\documentclass{amsart}

\usepackage{amsfonts}
\usepackage{amsmath}
\usepackage{amsthm}
\usepackage{amssymb}
\usepackage[english]{babel}
\usepackage{graphics}
\usepackage{latexsym}
\usepackage{longtable} \usepackage{mathrsfs}
\usepackage{graphicx}
\usepackage{wrapfig} %\usepackage{txfonts}
\usepackage[pdftex,colorlinks=true]{hyperref}

\newtheorem{theorem}{Theorem}
\newtheorem{corollary}[theorem]{Corollary}
\newtheorem{lemma}[theorem]{Lemma}

\newtheorem{claim}[theorem]{Claim}
\newtheorem{example}[theorem]{Example}

\theoremstyle{definition}
\newtheorem{definition}[theorem]{Definition}
\newtheorem{remark}[theorem]{Remark}

\newcommand{\mH}{\mathcal{H}}

\newcommand{\A}{\textbf{A}}
\newcommand{\R}{\mathbb{R}}

\newcommand{\mB}{\mathbb{B}}
\newcommand{\X}{\textbf{X}}

\newcommand{\D}{\mathscr{D}}

\newcommand{\noi}{\noindent}
\newcommand{\ms}{\medskip}

\newcommand{\al}{\alpha}
\newcommand{\be}{\beta}
\newcommand{\ga}{\gamma}
\newcommand{\de}{\delta}
\newcommand{\De}{\Delta}
\newcommand{\e}{\varepsilon}

\newcommand{\Si}{\Sigma}
\newcommand{\la}{\lambda}

\newcommand{\ka}{\kappa}
\newcommand{\Ga}{\Gamma}
\newcommand{\Om}{\Omega}

\newcommand{\lharpoonup}{-\!\!\!\!\rightharpoonup}

\newcommand{\larrow}{\longrightarrow}
\newcommand{\ot}{\otimes}

\newcommand{\lmapsto}{\longmapsto}
\newcommand{\ri}{\rightarrow}
\newcommand{\LL}{\text{\LARGE$\llcorner$}}
\newcommand{\p}{\partial}
\newcommand{\sub}{\subseteq}
\newcommand{\set}{\setminus}
\newcommand{\by}{\times}

\newcommand{\diam}{\textrm{diam}}

\newcommand{\spn}{\textrm{span}}

\newcommand{\bt}{\begin{theorem}}\newcommand{\et}{\end{theorem}}
\newcommand{\bd}{\begin{definition}}\newcommand{\ed}{\end{definition}}
\newcommand{\bl}{\begin{lemma}}\newcommand{\el}{\end{lemma}}
\newcommand{\beq}{\begin{equation}}\newcommand{\eeq}{\end{equation}}
\newcommand{\bc}{\begin{claim}}\newcommand{\ec}{\end{claim}}
\newcommand{\bex}{\begin{example}}\newcommand{\eex}{\end{example}}
\newcommand{\bcor}{\begin{corollary}}\newcommand{\ecor}{\end{corollary}}
\newcommand{\bp}{\begin{proof}}\newcommand{\ep}{\end{proof}}

\newcommand{\BPC}{\medskip \noindent \textbf{Proof of Claim} }

\newcommand{\BPT}{\medskip \noindent \textbf{Proof of Theorem} }

\numberwithin{equation}{section}
\begin{document}

\title[Linear Degenerate Elliptic Systems with Constant Coefficients]{On Linear Degenerate Elliptic PDE Systems with Constant Coefficients}

\author{Nikos Katzourakis}
\address{Department of Mathematics and Statistics, University of Reading, Whiteknights, PO Box 220, Reading RG6 6AX, Berkshire, UK}
%\curraddr{Department of Mathematics and Statistics, Case Western
%Reserve University, Cleveland, Ohio 43403}
\email{n.katzourakis@reading.ac.uk}

%    General info
%\subjclass[2010]{Primary 35J46, 35J47, 35J60; Secondary 35D30, 32A50, 32W50}

\date{}



\begin{abstract} Let $\A$ be a symmetric convex quadratic form on $\R^{Nn}$ and $\Om\sub \R^n$ a bounded convex domain. We consider the problem of existence of solutions $u: \Omega \subset \mathbb{R}^n \longrightarrow \mathbb{R}^N$ to the problem
\[ \tag{1} \label{1}
\begin{split}
\sum_{\be=1}^N\sum_{i,j=1}^n\A_{\al i \be j}\, D^2_{ij}u_\be \,=\,  f_\al, \text{ in }\Omega,\quad \  u\,=\, 0, \text{  on }\partial \Omega,
\end{split}
\]
when $ f\in L^2(\Omega,\R^N)$.  \eqref{1} is degenerate elliptic and it has not been considered before without the assumption of strict rank-one convexity. In general, it may not have even distributional solutions.  By introducing an extension of distributions adapted to \eqref{1}, we prove existence, partial regularity and by imposing an extra condition uniqueness as well. The satisfaction of the boundary condition is also an issue due to the low regularity of the solution. The motivation to study \eqref{1} and the method of the proof arose from recent work of the author \cite{K4} on generalised solutions for fully nonlinear systems.
\end{abstract}

\maketitle

\section{Introduction and the main result} \label{section1}

Let $n,N\geq 1$ be integers and consider a convex symmetric quadratic form $\A$ on the matrix space $\R^{Nn}$, that is
\beq \label{1.1}
\A_{\al i \be j}\, =\, \A_{\be j \al i},\ \ \ \ \sum_{\al,\be,i,j}\A_{\al i \be j }\, Q_{\al i}\, Q_{\be j}\, \geq\, 0, \quad  Q\in \R^{Nn}.
\eeq
Let also $\Om\Subset \R^n$ be a bounded convex domain. We will follow the convention that Greek indices $\al,\be,\ga,...$  run in $\{1,...,N\}$  and Latin indices $i,j,k,...$ run in $\{1,...,n\}$, even if their domain is not explicitly mentioned. In this note we consider the question of existence of solutions $u: \Omega \subset \R^n \larrow \R^N$ to the Dirichlet problem
\beq \label{1.2}
\left\{
\begin{split}
\sum_{\be,i,j} \A_{\al i \be j}\, D^2_{ij}u_\be\, &=\, f_\al, \ \text{ in }\Om, \quad \al=1,...,N,\\
u\,&=\, 0, \quad \text{ on }\p\Om, 
\end{split}
\right.
\eeq
when $ f\in L^2(\Omega,\R^N)$. The operators $D_iu_\al$, $D^2_{ij}u_\al$ will denote the $i$-th and $ij$-th partial derivatives of first and second order respectively of the $u_\al$ component of the map $u=(u_1,...,u_N)^\top$. 

It is well known that at least formally, \eqref{1.2} is the Euler-Lagrange equation of the convex functional
\beq  \label{1.3}
E(u,\Om)\, :=\, \int_\Om  \Bigg\{\frac{1}{2}\sum_{\al,\be,i,j} \A_{\al i \be j }\, D_iu_\al \, D_j u_\be  \ +\ \sum_{\ga} f_\ga\, u_\ga \Bigg\}
\eeq
in the Sobolev space $H^1_0(\Om,\R^N)$. If $\A$ is strictly rank-one convex on $\R^{Nn}$, that is when the following Legendre-Hadamard condition holds
\beq \label{1.4}
\sum_{\al,\be,i,j}\A_{\al i \be j }\, \eta_\al \, a_i\, \eta_\be\, a_j \, \geq\, c\, |\eta|^2|a|^2, \quad  \eta \in \R^N,\ a\in \R^n, 
\eeq
for some $c>0$ (where $|\cdot|$ denotes the Euclidean norm on both $\R^N$ and $\R^n$), then it is textbook material that the problem \eqref{1.2} has a unique weak solution which is a minimiser of \eqref{1.3} in the space $H^1_0(\Om,\R^N)$, see e.g.\ \cite{D, GM}. Moreover, by standard regularity results it follows that the solution is actually strong, lies in the space $(H^2\cap H^1_0)(\Om,\R^N)$ and satisfies the system a.e.\ on $\Om$ (e.g., \cite{GM}).

The primary advance in this paper is that we prove existence of solution to \eqref{1.2} without assuming the standard strict ellipticity condition \eqref{1.4}, but instead only the \emph{degenerate ellipticity} condition \eqref{1.1} and an extra constraint on the (nontrivial) nullspace of $\A$ which we explain later. Without \eqref{1.4} the functional \eqref{1.3} is convex but \emph{non-coercive} and standard variational/PDE methods fail. This is not a technical weakness, since as we show by examples (Ex. \ref{ex1}, \ref{ex2}) the solution in this case does not exist as an element of any standard Sobolev space and may not exist not even in the sense of distributions unless a compatibility condition is satisfied. To the best of our knowledge, the problem \eqref{1.2} has not been considered before without the assumption of strict rank-one convexity. 

The idea of the proof is based on the vanishing viscosity approximation of \eqref{1.2} by the strictly elliptic systems
\beq \label{1.5}
\left\{
\begin{split}
\sum_{\be,i,j} \big(\A_{\al i \be j}+\e \de_{\al \be}\de_{ij}\big)\, D^2_{ij}u^\e_\be\, &=\, f_\al, \ \text{ in }\Om, \quad \al=1,...,N,\ \e>0,\\
u^\e\,&=\, 0, \quad \text{ on }\p\Om, 
\end{split}
\right.
\eeq
and on the derivation of \emph{partial estimates along rank-one directions which are stable} as $\e\ri 0$.  By introducing an appropriate variant of the Distributional solutions adapted to the degeneracy of $\A$, we prove the existence of solution in this sense. We note that the satisfaction of the boundary condition is also a serious issue under the low regularity of the solution since the solution fails in general to be in $W^{1,1}_{\text{loc}}(\Om,\R^N)$ and there is no general trace operator for $L^1_{\text{loc}}(\Om,\R^N)$ mappings.

Before stating our existence result we need some preparation. Let $\A$ be a given tensor which satisfies \eqref{1.1} and will be fixed for the rest of the paper. The notation
\[
N\big(\A : \R^{Nn}\ri \R^{Nn} \big)
\]
will be used to denote the nullspace of $\A$ when $\A$ acts as a mapping
\[
\R^{Nn}\,\ni\, Q \lmapsto \A Q\, := \sum_{\al,\be,i,j}\big( \A_{\al i \be j }Q_{\be j}\big)\, e^\al \ot e^i\, \in \, \R^{Nn}.
\]
Evidently, $\{e^i\}$, $\{e^\al\}$ and $\{e^\al \ot e^i\}$ denote the Euclidean bases of $\R^n$, $\R^N$ and $\R^{Nn}$ respectively. Let us define the vector spaces
\beq \label{1.6}
\begin{split}
\Pi\, &:=\, N\big(\A : \R^{Nn}\ri \R^{Nn} \big)^\bot \hspace{30pt} \sub \R^{Nn},\\
\Si\, &:=\, \spn[\big\{\eta  \in \R^N\, :\, \eta \ot a \in \Pi \big\}] \  \sub \R^{N}.
\end{split}
\eeq
The space $\Pi$ is the orthogonal complement of the nullspace of $\A$ (namely, the range) and contains the ``rank-one directions of strict ellipticity", that is
\[
\sum_{\al,\be,i,j}\A_{\al i \be j }\, \eta_\al \, a_i\, \eta_\be\, a_j \, \geq\, c\, |\eta|^2|a|^2, \quad \ \eta \ot a \in \Pi\sub \R^{Nn}.
\]
We will follow the convention that the same letters $\Pi,\Si$ will denote the subspaces as well as the orthogonal projections on them. The meaning will be clear from the context, for example the projection map satisfies $\Si=\Si^\top=\Si^2\geq 0$ etc. The functional spaces of mappings $f:\Om \sub \R^n \larrow \Si\sub \R^N$ valued in $\Si$ will be denoted by $L^p(\Om,\Si)$, $W^{1,p}(\Om,\Si)$, $C^\infty(\Om,\Si)$, etc. Note also that we have $\Pi^\bot=\{0\}$ if and only if $\A$ defines a \emph{strictly} convex quadratic form, whence we also have $\Si=\R^N$ and $\Pi=\R^{Nn}$ in this case.  Let now
\[
\D(\Om,\Si) \, := \, C^\infty_c(\Om,\Si)
\]
be the space of ``test maps" valued in the subspace $\Si\sub \R^N$. We consider the space of Distributions ``valued in $\Si$", namely the dual space of $\D(\Om,\Si)$
\[
\D'(\Om,\Si)\, :=\, \big(\D(\Om,\Si)\big)^*.
\] 
We consider both spaces $\D, \D'$ as being equipped with their usual topologies (which we will not need, so we refer to \cite{F} for their definition).

\noi \textbf{Definition.} \emph{We will say that the map $u : {\Om}\sub \R^n \larrow \R^N$ is a \emph{Distributional solution in $\D'(\Om,\Si)$ of the PDE system
\[
\sum_{\be,i,j} \A_{\al i \be j}\, D^2_{ij}u_\be\, =\, f_\al, \ \text{ in }\Om, 
\]
when $u\in L^1_\text{loc} (\Om,\R^N)$ and for all $\phi \in \D(\Om,\Si)$ we have
\beq \label{1.7}
\int_\Om \sum_{\al,\be,i,j} \A_{\al i \be j}\, u_\be\, D^2_{ij}\phi_\al\, =\, \int_\Om \sum_\al f_\al \phi_\al.
\eeq}}

The following is our main result.

\begin{theorem}[Existence-Uniqueness-Partial regularity] \label{theorem} Let $n,N\geq1$ with $\Om\Subset \R^n$ a strictly convex bounded domain. Suppose that $\A$ is a quadratic form on $\R^{Nn}$ which satisfies \eqref{1.1}. We assume that the vector space $\Pi \sub \R^{Nn}$ (i.e.\ the orthogonal complement of the nullspace of $\A$, given by \eqref{1.6}) is spanned by rank-one directions. Then, for any $f\in L^2(\Om,\Si)$, the problem 
\[
\left\{
\begin{split}
\sum_{\be,i,j} \A_{\al i \be j}\, D^2_{ij}u_\be\, &=\, f_\al, \ \text{ in }\Om, \quad \al=1,...,N,\\
u\,&=\, 0, \quad \text{ on }\p\Om, 
\end{split}
\right.
\]
has a unique Distributional solution $u : \overline{\Om}\sub \R^n \larrow \R^N$ in $\D'(\Om,\Si)$ which lies in $L^2(\Om,\Si)$. In addition, $u$ is $\mH^{n-1}\LL \p U$-measurable on the boundary of any strictly convex subdomain $U\sub \Om $ and $u=0$ $\mH^{n-1}$-a.e.\ on $\p\Om$. 

Moreover, certain projections along rank-one directions of the Distributional gradient $Du$ exist as $L^2$ functions: for any $\eta \ot a \in \Pi$, the directional weak derivative of the projection $D_a(\eta \cdot u)$ exists in $L^2(\Om)$.
\end{theorem}

In the above statement, ``$\mH^{n-1}$" denotes the $(n-1)$-Hausdorff measure, ``$\cdot$" denotes the Euclidean inner product and ``$D_a$" is the standard directional derivate along $a$. We remark that $f$ must be valued in the subspace $\Si$ since this is a necessary \emph{compatibity condition} arising from the degenerate nature of the PDE. 

The following considerations show that the results and the assumptions of Theorem \ref{theorem} are optimal.

\begin{example}[Compatibity condition] \label{ex1} The compatibility condition in Theorem \ref{theorem} which requires that $f$ must be valued in $\Si$ is necessary: the degenerate $2\by2$ system
\[
\left\{
\begin{split}
\De u_1(x_1,x_2)\, &=\, f_1(x_1,x_2),\\
0\, &=\, f_2(x_1,x_2),
\end{split}
\right.
\]
has no solution whatsoever in any weak sense unless $f_2\equiv 0$.
\end{example}

\begin{example}[Partial regularity] \label{ex2} In general the solution we obtain in Theorem \ref{theorem} can not be a Sobolev function. Let $\Om=\mB_1(0)\sub \R^2$ be the unit disc centred at the origin and choose a function $f\in C^0(\overline{\Om})$ which is not weakly differentiable with respect to $x_1$ for any $x_2$. Then, the Dirichlet problem for the following degenerate elliptic single equation 
\[
\left\{
\begin{split}
D^2_{22}u(x_1,x_2)\, &=\, f(x_1,x_2), \text{ on }\Om,\\
u\,&=\, 0, \hspace{34pt} \text{ on }\p\Om,
\end{split}
\right.
\]
has the solution
\[
u(x_1,x_2)\, =\, -h(x_1,x_2)\ +\, \int_{-\infty}^{x_2}\int_{-\infty}^{t_2}f(x_1,s_2)\,ds_2\,dt_2.
\]
In the above, $h$ is the function 
\[
\begin{split}
h(x_1,x_2)\, :=\, &\left(\frac{g\big(x_1, \sqrt{1-x_1^2} \big) - g\big(x_1, -\sqrt{1-x_1^2} \big) } {2\sqrt{1-x_1^2}} \right)x_2\\
&\ \ \ \ \ \ \ \ +\, \frac{g\big(x_1, \sqrt{1-x_1^2} \big) + g\big(x_1, -\sqrt{1-x_1^2} \big) } {2} 
\end{split}
\]
where $g \in C^0(\p \Om)$ is the function which is given by 
\[
g(x_1,x_2)\, :=\, \int_{-\infty}^{x_2}\int_{-\infty}^{t_2}f(x_1,s_2)\,ds_2\,dt_2, \quad x_1^2+x_2^2=1.
\]
In view of our choice of $f$, the solution $u$ is not in $W^{1,1}_\text{loc}(\Om)$.
\end{example}

\begin{remark}[Nonuniqueness on the subspace of  ``degeneracies"] If we do not require the generalised solution we obtain in Theorem \ref{theorem} to satisfy $\Si^\bot u \equiv 0$ then it may not be unique (unless $\A$ is strictly elliptic, in which case we have $\Si^\bot =\{0\}$). Any extension of $u$ from $L^2(\Om,\Si)$ to $L^2(\Om,\R^N)$ which satisfies the boundary condition is also a solution. For instance,
\[
\tilde{u}\, :=\, \Si u\ +\ \Si^\bot g, \quad g\in C_0^0(\Om),
\]
is also a solution of the same Dirichlet problem for \emph{any} $g$. 
\end{remark}

The motivation to study the problem \eqref{1.2} and the method of proof come from the very recent paper of the author \cite{K4} and its companion paper \cite{K5}. In \cite{K4} we proposed a new duality-free theory of generalised solutions which applies to fully nonlinear PDE systems. This approach allows for nonlinearities of any order and with discontinuous coefficients whilst the only a priori regularity requirement of the solution is measurability. The standing idea of integration-by-parts which applies only to divergence systems is replaced by a probabilistic representation of derivatives which do not exist in the classical sense. Technically, this is done by utilising Young (parameterised) measures to describe the limiting behaviour of difference quotients over the compactification of the ``state space", that is the space wherein the derivatives are valued. 

Among other existence results in these papers, in \cite{K4} we proved existence of a so-called ``diffuse solution" to the problem
\beq  \label{1.8}
\left\{
\begin{split}
F(\cdot,D^2u)\, &=\, f, \ \text{ in }\Om,\\
u\,&=\, 0, \ \text{ on }\p\Om, 
\end{split}
\right.
\eeq
where $F : \Om\by \R_s^{Nn^2}\larrow \R^N$ is a Carath\'eodory mapping and $\R^{Nn^2}_s$ is the Euclidean space wherein the hessians $D^2u$ of smooth maps are valued. The result for  \eqref{1.8} extends previous work of the author for \eqref{1.8} but on the realm of strictly elliptic systems and of strong solutions (see \cite{K1, K2, K3}). The crucial assumption for existence is a degenerate ellipticity condition which roughly requires F to be ``close" to a linear degenerate system of the form we consider herein. This system has been solved in the ``diffuse" sense in \cite{K4} under assumptions stronger than those we consider herein, as a stepping stone in order to solve \eqref{1.8}. To aim of this paper is to show that under the present \emph{weaker assumptions}, the PDE system \eqref{1.2} has solutions in a certain distributional sense as well. 

A particular difficulty is the satisfaction of the boundary condition. In fact, the \emph{only} reason that strict convexity of $\Om$ is needed is for the satisfaction of the boundary condition. The strictness is meant in the sense that $\p\Om$ contains no non-trivial straight line segment. Although in general there is no trace operator because the solution may not be in any Sobolev space, yet it has ``differentiable rank-one projections". This means that the gradient $Du$ does not exist as a whole, but only certain projections of it exist along rank one lines of $\R^{Nn}$. Surprisingly, this suffices for a partial trace operator to exist.

The following condition which was introduced in \cite{K4} provides a sufficient condition about when the assumption of Theorem \ref{theorem} that $\Pi$ is spanned by rank-one directions is satisfied (see \eqref{1.6}). 

\ms

\noi \textbf{Structural Hypothesis (SH)} \textit{The tensor $\A$ satisfies (SH) if it can be decomposed as
\[
\A_{\al i \be j}\, =\, B^1_{\al \be}A^1_{ij}\, +\, \cdots\, +\, B^N_{\al \be}A^N_{ij}
\]
and also: }

\noi \textit{a) The symmetric matrices $\{A^1,...,A^N\} \sub \R^{n^2}$ are non-negative and the eigenspaces
\[
N\Big(A ^\ga-\la^\ga_+ I : \R^n\ri \R^n\Big) 
\]
intersect for all $\ga=1,...,N$ along a common line. Here $\la^\ga_+$ denotes the smallest positive eigenvalue of $A^\ga$.}

\noi \textit{b) The symmetric matrices $\{B^1,...,B^N\} \sub \R^{N^2}$ are non-negative and have mutually orthogonal ranges.}

\ms

Note that (SH)  trivialises when either $N=1$ or $n=1$ since any symmetric non-negative matrix satisfies it. Although (SH) is quite restrictive, by its constructive nature is is evident how to demonstrate nontrivial examples when $n,N\geq 2$.

\section{Proof of the main result} \label{section2}

\BPT \ref{theorem}. The proof is based on the approximation by the strictly elliptic systems \eqref{1.5} as $\e\ri0$. 

\ms

\textbf{Step 1}. To begin with, fix $f\in L^2(\Om,\Si)$, $\e>0$ and consider the functional
\beq \label{2.1}
E_\e(u,\Om)\, :=\, \int_\Om \Bigg\{ \frac{1}{2}\sum_{\al,\be,i,j} \A_{\al i \be j }\, D_iu_\al \, D_j u_\be \ +\,\frac{\e}{2}|Du|^2\, +  \sum_{\ga} f_\ga\, u_\ga \Bigg\}
\eeq
placed in $H^1_0(\Om,\R^N)$. Standard semicontinuity and regularity results (see e.g.\ \cite{D, GM}) imply that the problem \eqref{1.5} has a unique strong solution $u^\e \in (H^2\cap H^1_0(\Om,\R^N)$ which also is a global minimiser of \eqref{2.1} over $H^1_0(\Om,\R^N)$. 

\ms

\textbf{Step 2.} Now we obtain estimates stable in $\e$. 

\ms

We begin with some algebraic consequences of our assumptions. By \eqref{1.6}, \eqref{1.1} and the Spectral theorem applied to the symmetric linear map $\A :\R^{Nn}\larrow \R^{Nn}$, we have that there exists $\nu>0$ such that
\begin{align}
\Pi \A \Pi \, =\, \A\, &=\, \A \Pi\,=\, \Pi\A, \label{2.2}\\
\sum_{\al ,i,\be,j}\A_{\al i \be j}\, Q_{\al i}\, Q_{\be j}\, &\geq\, \nu\, |\Pi Q|^2, \quad Q\in \R^{Nn}.  \label{2.3}
\end{align}
In \eqref{2.3}, $|\cdot|$ denotes the Euclidean norm on $\R^{Nn}$, namely $|Q|^2=\sum_{\al ,i}Q_{\al i}Q_{\al i}$. The identities \eqref{2.2} say that $\A$ commutes with the projection map $\Pi$ on its range. \eqref{2.3} is a consequence of the fact that the restriction $\A\big|_{\Pi}$ on the range is invertible. Moreover, we have that
\beq \label{2.4}
\X \in \R^{Nn^2} \ \ \Longrightarrow \ \ \sum_{\al,\be,i,j} \big(\A_{\al i \be j}\,  \X_{\be i j}\big)e^\al \ \in \Si \sub \R^{Nn}.
\eeq
\eqref{2.4} says that the image of the linear map $\A : \R^{Nn^2}\larrow \R^N$ is a subspace of $\Si$. This is a result of the degeneracy of $\A$ and can be seen as folllows: every $\X \in \R^{Nn^2}$ can be written as
\[
\X\, =\sum_{A=1}^{Nn^2} \xi^A \ot d^A \ot b^A,\ \quad \xi^A\in \R^N,\  \, d^A,b^A\in \R^n. 
\]
Then, by \eqref{2.2} we have (if the components of the projection $\Pi$ are denoted by $\Pi_{\al i \be j}$ and we use the symmetry of $\Pi$)
\[
\begin{split}
\eta_\al \, :=\, \sum_{\be,i,j}  \A_{\al i \be j}\,  \X_{\be i j} \, &=\, \sum_{A=1}^{Nn^2} \sum_{\be,i,j}  \A_{\al i \be j}\, \xi_{\be}^A \,d^A_i \, b_j^A  \\
&=\, \sum_{ \ka,k,i}\left\{ \sum_{A=1}^{Nn^2} \sum_{\be,j}    \A_{\ka k \be j} \,    \, \xi_{\be}^A \, b^A_j  \right\}  \Pi_{\ka k \al i}\, d_i^A .
\end{split}
\]
By the definition of $\Si$ in \eqref{1.6},  the above says that the vector $\eta :=\sum_\al \eta_\al e^\al$ belongs to $\Si$ and hence \eqref{2.4} follows.

Now, since $u^\e$ is a  minimiser, by \eqref{2.2} and by the compatibility condition which says that $\Si f =f$ (because $f$ is valued in $\Si\sub \R^N$), we have
\[
\begin{split}
0\, =\, E_\e(0,\Om) \, &\geq\, E_\e(u^\e,\Om)\\
&\geq \, \frac{1}{2}\int_\Om \sum_{\al,\be,i,j} \A_{\al i \be j }\, D_iu^\e_\al \, D_j u^\e_\be \  -\,  \int_\Om \sum_{\ga} (\Si f)_\ga\, u^\e_\ga .
\end{split}
\]
By \eqref{2.3} and by the symmetry of the projection map $\Si$,  for any $\de>0$ small we deduce 
\[
\begin{split}
0\, &\geq \ \frac{1}{2}\int_\Om \sum_{\al,\be,i,j} \A_{\al i \be j }\, D_iu^\e_\al \, D_j u^\e_\be \  -\,  \int_\Om \sum_{\ga} (\Si f)_\ga\, u^\e_\ga  \\
&\geq \   \frac{\nu}{2} \int_\Om \big|\Pi Du^\e\big|^2  \, -\,  \int_\Om \sum_{\ga}f_\ga\,  (\Si u^\e)_\ga 
\end{split}
\]
which gives
\beq \label{2.5}
\begin{split}
0\,\geq \  \frac{\nu}{2}\int_\Om \big|\Pi Du^\e\big|^2 \, -\,  \frac{4}{\de}\int_\Om |f|^2\, -\, \de\int_{\Om} |\Si u^\e|^2. \\
\end{split}
\eeq
Now we need a generalisation of the Poincar\'e inequality and of the trace operator. The rest of the proof follows similar line to those of \cite{K4}, by we provide all the details for the convenience of the reader.

\begin{claim}[Partial Poincar\'e inequality, cf. \cite{K4}] \label{Poincare} Let $\Om,\Pi,\Si$ be as in Theorem \ref{theorem}. Then, there exists $C=C(\Om,n,N)>0$ depending only on the diameter of $\Om$ and on the dimensions $n,N$ such that, for any $u\in H^1_0(\Om,\R^N)$ we have the estimate
\[
\|\Si  u\|_{L^2(\Om)}\, \leq\, C\, \big\|\Pi Du \big\|_{L^2(\Om)}.
\]
\end{claim}

We note that Claim \ref{Poincare} is actually true for any bounded open domain $\Om$.

\BPC \ref{Poincare}. Fix vectors $e\in \R^N$ and $\eta \in \R^N$ and let us denote by $e^\bot$ the hyperplane normal to $e$. For any $y\in e^\bot$, we set  (see Figure 1)
\[
I^{y,e}\,:=\, \big\{t\in \R\, \big|\, y+te\in \Om\big\}, \ \ \ \ \Om^e\,:=\, \Big\{y\in e^\bot\, \big|\ \exists\,t\in\R\, :\, y+te\in \Om\Big\}.
\]
We fix a function $u\in C^1_0(\Om,\R^N)$ and some $x= y+te\in \Om$. Then, for the projection $\eta \cdot u$ along $\eta$ we have
\[
\big| (\eta \cdot u)(y+te)\big|^2\, \leq\, |I^{y,e}|\int_{I^{y,e}} \big| e\cdot D(\eta \cdot u)(y+\la e)\big|^2 d\la.
\]
By integration with respect to $t\in I^{y,e}$ and $y\in \Om^e$, Fubini's theorem implies
\[
\begin{split}
\int_\Om \big| (\eta \cdot u)(x)\big|^2dx\, &\leq\, \int_{\Om^e}\left(|I^{y,e}|^{2}\int_{I^{y,e}} \big| e\cdot D(\eta \cdot u)(y+\la e)\big|^2 d\la\right)d\mH^{n-1}(y)\\
&\leq\, \diam(\Om)^2\int_{\Om}  \big|\eta \ot e: Du(x)\big|^2 dx.
\end{split}
\]
\[
\underset{ \text{Figure 1.} }{\includegraphics[scale=0.16]{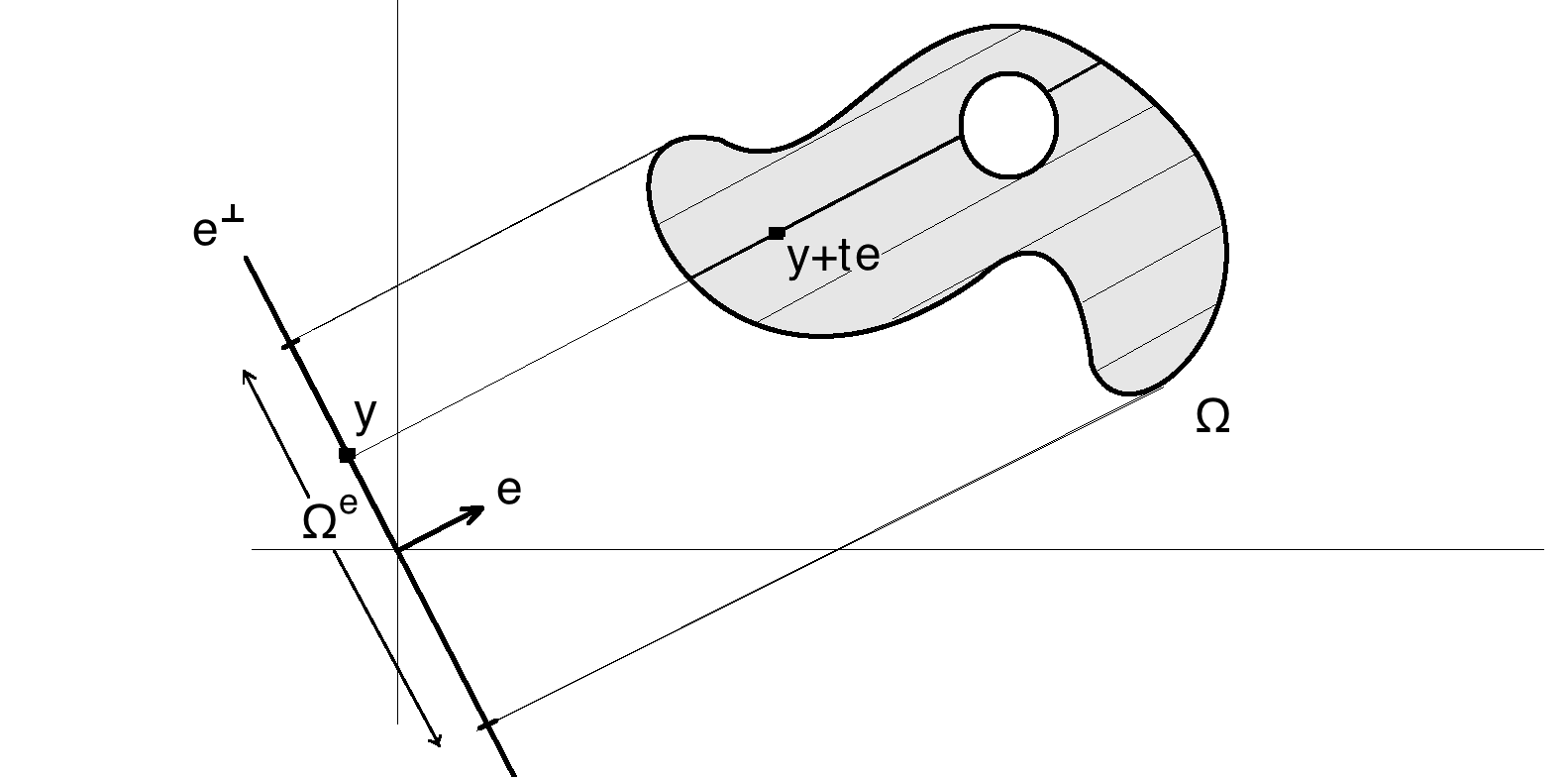}}
\]
Since by assumption $\Pi \sub \R^{Nn}$ is spanned by rank-one directions of the form $\eta \ot e$, by the definition of $\Si$ in \eqref{1.6} the desired estimate follows by considering a basis of $\Pi$ consisting of rank-one matrices and by using a standard density argument in the Sobolev norm.             \qed

\ms

Finally, by using Claim \ref{Poincare} and \eqref{2.5}, by choosing $\de>0$ small enough we have the uniform in $\e$ estimate
\beq \label{2.6}
\|\Si  u^\e\|_{L^2(\Om)}\, + \, \big\|\Pi Du^\e \big\|_{L^2(\Om)}\, \leq\, C\| f\|_{L^2(\Om)}
\eeq
for all $\e>0$. By \eqref{2.6} and by weak compactness, there exist maps $u \in L^2(\Om,\Si)$ and $U\in L^2(\Om,\Pi)$ such that
\beq  \label{2.6a}
\begin{split}
\Si  u^\e &\lharpoonup u,\ \text{ in }L^2(\Om,\Si),\\
\Pi  Du^\e &\lharpoonup U,\ \text{ in }L^2(\Om,\Pi),\\
\end{split}
\eeq
along a sequence $\e=\e_k\ri 0$ as $k\ri \infty$.

\ms

\textbf{Step 3.} Now we prove existence of a solution to \eqref{1.2}. Since for any $\e>0$ $u^\e$ is a strong solution of \eqref{1.5}, we have 
\[
\sum_{\be,i,j} \big(\A_{\al i \be j}+\e \de_{\al \be}\de_{ij}\big)\, D^2_{ij}u^\e_\be\, =\, f_\al, 
\]
a.e.\ on $\Om$ for all indices $\al$. Moreover, $u^\e=0$ $\mH^{n-1}$-a.e.\ on $\p \Om$ by standard results on the pointwise properties of Sobolev functions and the regularity of the domain (see e.g.\ \cite{EG}). Let us fix a test map $\phi \in \D(\Om,\Si)$. Integration by parts gives
\beq \label{2.7}
\int_\Om \sum_{\al,\be,i,j} \A_{\al i \be j}\, u^\e_\be\, D^2_{ij}\phi_\al \, =\, \int_\Om \sum_\ga f_\ga \phi_\ga\, -\,\e  \int_\Om \sum_\la u^\e_\la\, \De\phi_\la.
\eeq
By \eqref{2.4} we have that 
\beq \label{2.8}
\sum_{ \be,i,j}  \A_{\al i \be j}\,  D^2_{i j}u^\e_\be   \, =\, \sum_{\la,\be,i,j} \big(\A_{\al i \be j}\,  D^2_{i j}u^\e_\be \big)\Si_{\al \la}
\eeq
and hence by using that $\Si \phi =\phi$ and \eqref{2.8}, \eqref{2.7} gives
\[
\int_\Om \sum_{\ka,\al,\be,i,j} \left(\A_{\al i \be j}\, D^2_{ij}\phi_\al \right)\big(\Si_{\be \ka} u^\e_\ka\big) \, =\, \int_\Om \sum_\ga f_\ga \phi_\al\, -\,\e  \int_\Om \sum_\la \big(\Si_{\la \ka}u^\e_\ka\big)\, \De\phi_\la.
\]
By letting $\e=\e_k\ri 0$, the convergences of \eqref{2.6a} imply that the previously obtained limit map $u \in L^2(\Om,\Si)$ is a Distributional solution in $\D'(\Om,\Si)$ of the system
\[
\sum_{ \be,i,j}  \A_{\al i \be j}\,  D^2_{i j}u_\be   \, =\, f_\al.
\]

\textbf{Step 4.} We now consider the problem of the satisfaction of the boundary condition. The next result shows that for strictly convex domains we have a partial trace operator under the assumptions of Theorem \ref{theorem}.

\begin{claim}[Trace operator, cf. \cite{K4}] \label{trace} Let $\Om,\Pi,\Si$ be as in Theorem \ref{theorem}. Then, there exists a closed $\mH^{n-1}$-nullset $E\sub \p\Om$ such that, for any $\Ga\Subset \p\Om\set E$, we can find $C=C(n,\Ga)>0$:
\[
\| \Si u \|_{L^2(\Ga,\mH^{n-1})} \, \leq\, C \,\Big( \| \Si u  \|_{L^2(\Om)}  \,+\, \big\| \Pi Du \big\|_{L^2(\Om)} \Big),
\]
for all maps $u\in H^1(\Om,\R^N)$. In addition, we have that
\[
E\,=\, \Big\{ x\in \p\Om\ :\ \big( x+\spn[e]\big)\cap \p\Om = \{x\} \Big\}.
\]
\end{claim}

We note that Claim \ref{trace} is a minor extension of standard results (see e.e.\ \cite{E, EG}).

\BPC \ref{trace}. Suppose $E\sub \p\Om$ is the closed set defined in the statement of the claim. By the strict convexity of $\Om$, it can be seen that $\mH^{n-1}(E)=0$. Let us fix a function $u\in C^1(\overline{\Om},\R^N)$ and a unit rank-one matrix $\eta \ot e\in \Pi\sub \R^{Nn}$.  We cover $\p\Om\set E$ by a sequence of open cubes $\{Q_j\}_1^\infty$ whose sides are orientated parallel to $\{e,e^\bot\}$ (see Figure 2). For every cube $Q_j$, we consider the sets 
\[
\Om_j\, :=\, Q_j \cap \Om, \quad \Ga_j \, :=\, Q_j \cap \p\Om. 
\]
Let us fix a triplet $(Q_j,\Om_j, \Ga_j)$ and assume that $e$ points towards the interior of $\Om_j$ (for otherwise we may replace it by $-e$). We may also restrict $\eta \cdot u ,D(\eta \cdot u)$ on $\Om_j$ and define 
\[
(\eta \cdot u)_j \, :=\, (\eta \cdot u)\chi_{\Om_j}, \quad (D(\eta \cdot u))_j\, :=\, D(\eta \cdot u)\chi_{\Om_j}. 
\]
Then, standard estimates imply that for any $x\in \Ga_j$  we have
\[
\big|(\eta \cdot u)(x)\big|^2\, \leq\, C\left( \int_0^\infty \big| (\eta \cdot u)_j(x+te) \big|^2dt\, +\,  \int_0^\infty \big|  \big(D(\eta \cdot u)\big)_j(x+te)\cdot e \big|^2dt\right).
\]
\[
\underset{ \text{Figure 2.} }{\includegraphics[scale=0.18]{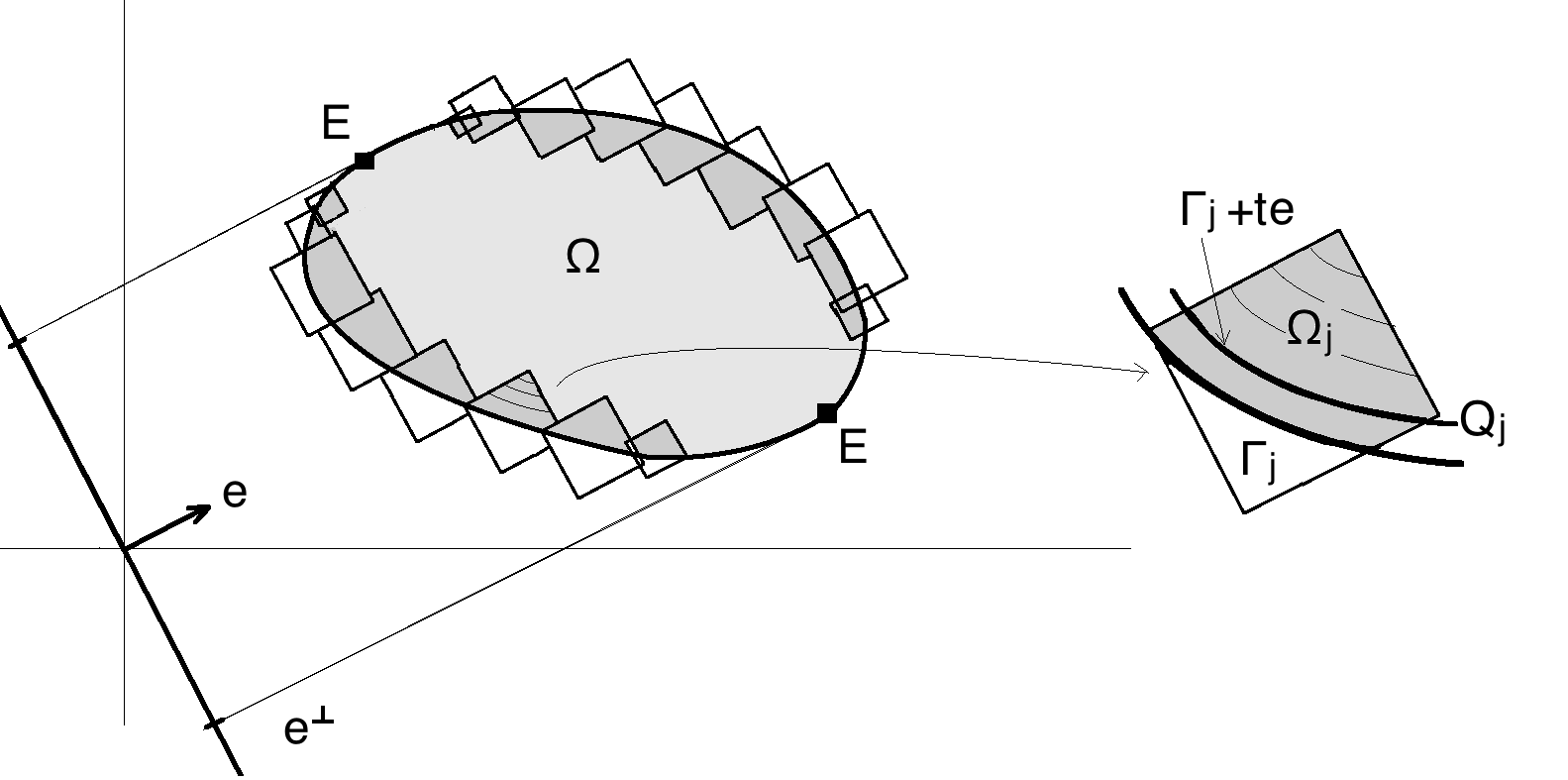}}
\]
Suppose now that $F_j \in W^{1,\infty}(\R^n)$ is a Lipschitz function such that when restricted to $\Om_j$, it satisfies the following: for each $t\geq 0$, the level set $\{F_j=t\}$ coincides with the translate of a portion of the boundary $ \p\Om+te$.  Such a function $F_j$ can be given by
\[
F_j(z) \, :=\, \sup\{t>0: z\in ({\Om}+te)\cap\Om_j\},  \quad z\in \Om_j,
\]
and can be trivially extended to a Lipschitz function on $\R^n$. Then, we integrate over $x\in \Ga_j$ and apply the co-area formula (see e.g.\ \cite{EG}) and Fubini's theorem:
\[
\begin{split}
\int_{\Ga_j} \big|(\eta \cdot u)\big|^2 \, d  \mH^{n-1}& 
\leq\, C \int_\R  \int_{\Ga_j +te}  \Big( \big| (\eta \cdot u)_j(y) \big|^2+ \big| \big(D(\eta \cdot u)\big)_j(y)\cdot e \big|^2\Big) \,d\mH^{n-1}(y)\, dt\\
=&\ \, C  \int_\R  \int_{\{F_j=t\}}  \Big( \big| (\eta \cdot u)_j(y) \big|^2+  \big| \big(D(\eta \cdot u)\big)_j(y)\cdot e \big|^2\Big) \,d\mH^{n-1}(y)\, dt\\
=& \ \, C    \int_{\R^n} \big|DF_j(x)\big| \Big(\big| (\eta \cdot u)_j(x) \big|^2 +  \big| \big(D(\eta \cdot u)\big)_j(x)\cdot e\big|^2\Big) \,dx\\
\leq &\ \, C   \, \|DF_j \|_{L^\infty(\R^n)}   \int_{\Om_j}\Big(\big| (\eta \cdot u)(x) \big|^2+\big| \eta \ot e : Du(x)  \big|^2\Big) \,dx.
\end{split}
\]
By using the assumption that $\Pi$ has a basis consisting of rank-one matrices $\eta \ot e$ and  $\Si$ is spanned by the respective directions $\eta \in \R^N$, the rest of the proof is an obvious application of a standard argument of partitions of unity.      \qed

\ms

An application of Claim \ref{trace} shows that the Distributional solution $u$ is $\mH^{n-1}\LL \p\Om$-measurable and $u=0$ $\mH^{n-1}$-a.e.\ on $\p \Om$. In addition, it is $\mH^{n-1}$-measurable on the boundary of any strictly convex subdomain of $\Om$.  

\ms

\textbf{Step 5.} We now show the uniqueness of the solution. Suppose the problem has 2 solutions $u,v$ and set $w:=u-v$. Then, $w\in L^2(\Om,\Si)$ is a solution in $\D'(\Om,\Si)$
of 
\[ 
\begin{split}
\sum_{\be=1}^N\sum_{i,j=1}^n\A_{\al i \be j}\, D^2_{ij}w_\be \,=\,  0, \text{ in }\Omega,\quad \  w\,=\, 0, \text{  on }\partial \Omega.
\end{split}
\]
(The satisfaction of the boundary condition is considered in the sense of the statement of the theorem.) Then, $w$ can be extended on $\R^n\set \Om$ by zero as a function in $L^2(\R^n,\Si)$. By mollifying in the standard way, for any $\e>0$ the mollified solution $w^\e=w *\eta^\e$ satisfies 
\[
\sum_{\be=1}^N\sum_{i,j=1}^n\A_{\al i \be j}\, D^2_{ij}(w^\e)_\be \,=\,  0, \quad \text{ on }\R^n,
\]
and $w^\e \in C^\infty_c(\R^n,\Si)$. Integration by parts and application of the inequality \eqref{2.3} give
\[
\big\|\Si Dw^\e \big\|^2_{L^2(\R^n)} \, \leq\, \frac{1}{\nu}\int_{\R^n} \sum_{\al, i, \be, j} \A_{\al i \be j}\, D_i (w^\e)_\al\, D_j (w^\e)_\be\, =\, 0.
\]
FInally, by Claim \ref{Poincare} and by the compactness of the support of $w^\e$ we obtain $\Si w^\e=0$ on $\R^n$. Since $\Si^\bot w\equiv 0$, by letting $\e\ri 0$ we get $w\equiv 0$ and as such we infer that the solution of the problem is unique.

\ms

\textbf{Step 6.} In order to conclude it remains to show the weak differentiability of the projection $\eta \cdot u$ along the $a$-direction of $\R^n$ when $\eta \ot a \in \Pi$. This is a consequence of the convergences in \eqref{2.6a}: for any $\phi \in C^\infty_c(\Om)$, we have
\[
\begin{split}
\int_\Om (\eta \cdot u )D_a\phi\, & =\, \lim_{k\ri \infty}\int_\Om (\eta \cdot u^{\e_k} )D_a\phi\\
& =\, -\lim_{k\ri \infty}\int_\Om \phi\, D_a(\eta \cdot u^{\e_k}) \\
& =\, -\lim_{k\ri \infty}\int_\Om  \phi\, \big(\eta \ot a :D u^{\e_k}\big)
\end{split}
\]
and hence
\[
\begin{split}
\int_\Om (\eta \cdot u )D_a\phi\, =\, -\int_\Om \phi \, \big(\eta \ot a :U\big),
\end{split}
\]
for any $\phi \in C^\infty_c(\Om)$. As a consequence, we have that $D_a (\eta \cdot u )=\eta \ot a :U$ and also that $\eta \ot a :U$ is in $L^2(\Om)$. The proof of the theorem is completed.       \qed

\ms

\noi \textbf{Acknowledgement.} The author is indebted to the anonymous referee for the careful yet quick reading of this paper and for their suggestions. In particular, we would like to thank them for a correction in Example 3.

\ms

\bibliographystyle{amsplain}

\end{document}